\documentclass[12pt]{article}
 \usepackage{amsfonts}
\usepackage{amsmath,amssymb}
\usepackage{amscd}
\usepackage{amsthm}
\title{On hom-algebras with surjective twisting}
\author {Aron Gohr\\
{\small University of Luxembourg}\\
{\small Mathematics Research Unit}\\
{\small  162A, avenue de la faiencerie}\\
{\small  L-1511, Luxembourg}\\
 {\small  Grand-Duchy of Luxembourg.}\\
}

\newtheorem{definition}{Definition}
\newtheorem{proposition}{Proposition}
\newtheorem{lemma}{Lemma}

\newtheorem{remark}{Remark}
\begin {document}
\date{\today}

\maketitle
\begin{abstract}
A hom-associative structure is a set $A$ together with a binary operation $\star$ and a selfmap $\alpha$ such that an $\alpha$-twisted version of associativity is fulfilled. In this paper, we assume that $\alpha$ is surjective. We show that in this case, under surprisingly weak additional conditions on the multiplication, the binary operation is a twisted version of an associative operation. As an application, an earlier result \cite{FG2009} on weakly unital hom-algebras is recovered with a different proof. In the second section, consequences for the deformation theory of hom-algebras with surjective twisting map are discussed.
\end{abstract}
\section*{Introduction}
The study of \emph{hom-algebras} originates with \cite{HLS}, who introduced a notion of hom-Lie algebra in the context of deformation theory of Witt and Virasoro algebras. Later, this notion was generalized and transferred to other categories in \cite{MS2008}. Deformation theory of hom-associative algebras was first explored in \cite{MS2007}, where a Gerstenhaber-type notion of formal deformation for hom-associative algebras is introduced and the beginnings of a cohomology theory appropriate for studying such deformations are developed. In \cite{Yau:Homology}, using a different notion of deformation, new examples of hom-associative and hom-Lie algebras were constructed from associative respectively Lie algebras. In \cite{FG2009}, it was shown that in the case of \emph{unital} hom-associative algebras relatively innocent-looking conditions on the twisting map can force a hom-associative algebra to be associative. For instance, unital hom-associative algebras with surjective twisting are associative. Also a weakened notion of unitality for hom-associative algebras was investigated and it was found that under the assumption of weak unitality and bijective twisting, a hom-associative algebra, while not being necessarily associative, can always be constructed from an associative algebra by a generalization of one of the construction procedures in \cite{Yau:Homology}.\\
This paper is divided in two sections after this introduction. The first section extends and improves upon some of the findings of \cite{FG2009}. We prove a very general result (Proposition \ref{twisting-theorem}) about hom-associative structures with surjective twisting, which says essentially that either the multiplication on such a structure is in some way degenerate or it can be constructed from an associative structure as in \cite{Yau:Homology}. The notion of nondegeneracy of an algebraic structure used here will be made precise in the first section. The previous result on weakly unital hom-structures with bijective twisting from \cite{FG2009} is obtained as a special case. A second theorem is subsequently proven which shows that the assumption of \emph{surjective} instead of \emph{bijective} twisting map in Proposition \ref{twisting-theorem} is no real advantage over what was obtained in \cite{FG2009}, because in our situation the twisting map can be shown to be bijective anyway. However, the replacement of weak unitality with nondegeneracy as well as the different method of proof employed provide real progress over \cite{FG2009}.\\
The second section is devoted to a treatment of hom-deformation theory in the sense of \cite{MS2007}, in the special case where the hom-associative algebra to be deformed had a surjective twisting map and a nondegenerate multiplication. The key observation here is that both nondegeneracy and surjectivity of the twisting map are preserved under hom-associative deformation of hom-associative algebras. We use this observation to relate hom-associative deformations of an algebra $B$ which arises from a surjective twisting of an associative algebra $A$ to associative deformations of $A$.\\
We use similar conventions and notations as in \cite{FG2009}. Specifically, $k$ will always denote a commutative ring with unit, $K$ will be a field. 

\section{Hom-structures with surjective twisting}
A hom-associative structure is a set $A$ together with a multiplication $\star:A \times A \rightarrow A$ and a self-map $\alpha:A \rightarrow A$ such that the condition
\begin{displaymath}
\alpha(x) \star (y \star z) = (x \star y) \star \alpha(z)
\end{displaymath}
is fulfilled. Depending on the category under consideration, $\alpha$ is expected to satisfy other conditions as well. In general, the philosophy is that $\alpha$ should be a homomorphism for all functions and relations on our algebraic structure, except possibly for the multiplication $\star$. For example, in the case of hom-rings $\alpha$ is supposed an abelian group endomorphism, in the case of hom-$k$-algebras it is linear over the commutative base ring $k$. However, in the present section we need no such additional structures on $A$.\\
If we consider $A$ and $\alpha$ fixed, it is clear that the extent to which such an associativity condition restricts the possible choices of $\star$ is highly dependent on the choice of $\alpha$. In the case of hom-algebras for instance, it is possible to choose $\alpha = 0$ and obtain a hom-associative structure with every bilinear $\star:A \times A \rightarrow A$. On the other hand, if $\alpha = id$, we obtain the usual notion of associative algebras. In this section, we will study some aspects of what happens when $\alpha$ is assumed to be \emph{surjective}.
\subsection{Preliminaries and definitions}
We recall in the following the definition of a \emph{twist} from \cite{FG2009}. We also introduce the notion of an \emph{untwist}:
\begin{definition}
Let $(A, \star, \alpha)$ be a hom-associative structure. Then $A$ is called a \emph{twist} if there is an associative multiplication $\cdot: A \times A \rightarrow A$ such that $(A, \star, \alpha)$ arises from $(A, \cdot)$ by setting $x \star y := \alpha(x \cdot y)$. The structure $(A, \cdot)$ is called an \emph{untwist} of $(A, \star, \alpha)$.
\end{definition}
Note that any such multiplication must due to the hom-associativity of $A$ satisfy the hom-associativity-like condition
\begin{displaymath}
\alpha(\alpha(x) \cdot \alpha(y \cdot z)) = \alpha(\alpha(x \cdot y) \cdot \alpha(z)).
\end{displaymath}
If $\alpha$ is bijective and if there is an untwist of of $A$, then this untwist is obviously uniquely defined. We will therefore in this case talk of \emph{the} untwist of $(A, \star, \alpha)$.\\
In \cite{FG2009} it was shown that weakly unital hom-associative structures with bijective twisting are always twists, and indeed twists of unital structures. Equivalently, one could say that any unital structure $(A, \cdot, \alpha, 1)$ such that $\alpha$ is bijective and satisfies the hom-associativity-like relation
\begin{displaymath}
\alpha(\alpha(x) \cdot \alpha(y \cdot z)) = \alpha(\alpha(x \cdot y) \cdot \alpha(z)).
\end{displaymath}
is associative.\\
In the sequel, our goal is to show that the condition of weak unitality is not essential here.

\subsection{Surjective twistings: the main lemma}
Our first goal is to prove a technical lemma about hom-associative structures with surjective twisting. The motivation behind the introduction of this lemma is as follows: suppose that $(A, \star, \alpha)$ is a hom-associative structure with $\alpha$ surjective. Then ideally, we would like to be able to write the multiplication $\star$ as a twisting of an associative multiplication by $\alpha$, i.e. in the form $$x \star y := \alpha(x \cdot y)$$, where $\cdot:A \times A \rightarrow A$ is associative. Since by surjectivity of $\alpha$ there exist $\beta:A \rightarrow A$ with $\alpha \circ \beta = id_A$, a natural ansatz is to simply set $x \cdot y := \beta(x \star y)$ with such a $\beta$. The associativity condition $(x \cdot y) \cdot z = x \cdot (y \cdot z)$ is then the same as $x \star \beta(y \star z) = \beta(x \star y) \star z$, since $\beta$ is necessarily injective.\\
In general, $\beta$ cannot be chosen such that this associativity condition is fulfilled. However, the following weaker statement can be shown:
\begin{lemma} \label{haupt-hilfssatz}
Let $(A, \star)$ be a hom-associative structure with $\alpha$ surjective and let $\beta:A \rightarrow A$ be a map with $\alpha \circ \beta = id_A$. Then, the following associativity conditions are satisfied for all $a,b,x,y,z \in A$:
\begin{eqnarray}
a \star (b \star (x \star \beta(y \star z))) & = & a \star (b \star (\beta(x \star y) \star z)), \label{ass-1}\\
a \star ((x \star \beta(y \star z)) \star b) & = & a \star ((\beta(x \star y) \star z) \star b), \label{ass-2}\\
((x \star \beta(y \star z)) \star b) \star a & = & ((\beta(x \star y) \star z) \star b) \star a, \label{ass-3}\\
(b \star (x \star \beta(y \star z))) \star a & = & (b \star (\beta(x \star y) \star z)) \star a. \label{ass-4}
\end{eqnarray}
\end{lemma}
\begin{proof}
Equations \ref{ass-3} and \ref{ass-4} can be obtained from (Eq. \ref{ass-1}, \ref{ass-2}) by passing to the opposite hom-structure, i.e. to the hom-associative structure with multiplication $a \star^{op} b := b \star a$. \\
We first remark that we have
\begin{equation}
(\beta(x) \star y) \star z = x \star (y \star \beta(z)) \label{hilfsgleichung1}
\end{equation}
as was already shown in \cite{FG2009}. We also have for any $x,y,z,u \in A$ the identity
\begin{equation}
\alpha^2(x) \star ((y \star z) \star u) = \alpha(x \star y) \star (\alpha(z) \star u) \label{hilfsgleichung2}
\end{equation}
because of
\begin{eqnarray*}
\alpha^2(x) \star ((y \star z) \star u) & = & (\alpha(x) \star (y \star z)) \star \alpha(u)\\
& = & ((x \star y) \star \alpha(z)) \star \alpha(u) \\
& = & \alpha(x \star y) \star (\alpha(z) \star u).
\end{eqnarray*}
With this said, we are ready for a proof of (Eq. \ref{ass-1}):
\begin{eqnarray*}
a \star (b \star (x \star \beta(y \star z))) & \overset{\text{Eq. \ref{hilfsgleichung1}}}{=} & a \star ((\beta(b) \star x) \star (y \star z)) \\
& \overset{\alpha \circ \beta = id_A}{=} & \alpha^2(\beta^2(a)) \star ((\beta(b) \star x) \star (y \star z)) \\
& \overset{\text{Eq. \ref{hilfsgleichung2}}}{=} &  \alpha(\beta^2(a) \star \beta(b)) \star ((x \star y) \star \alpha(z)) \\
& \overset{\text{hom-ass.}}{=} & ((\beta^2(a) \star \beta(b)) \star (x \star y)) \star \alpha^2(z)\\
& \overset{\text{Eq. \ref{hilfsgleichung1}}}{=} & (\beta(a) \star (\beta(b) \star \beta(x \star y))) \star \alpha^2(z) \\
& \overset{\alpha \circ \beta = id_A}{=} & a \star ((\beta(b) \star \beta(x \star y)) \star \alpha(z)) \\
& \overset{\alpha \circ \beta = id_A}{=} & a \star (b \star (\beta(x \star y) \star z)).
\end{eqnarray*}
The proof of (Eq. \ref{ass-2}) follows the same method, starting from 
\begin{displaymath}
a \star ((\beta(x \star y) \star z) \star b) = a \star ((x \star y) \star (z \star \beta(b))). 
\end{displaymath}
The same trick as above of replacing $a$ by $\alpha^2(\beta^2(a))$ and using (Eq. \ref{hilfsgleichung2}) is applied to obtain a sub-term of the form $\alpha(y) \star (z \star \beta(b))$, which simplifies to $(y \star z) \star b$. The same simplifying steps as above then yield (Eq. \ref{ass-2}).
\end{proof}
\subsection{Nondegenerate multiplications}
Our first main result will be that a hom-associative structure $(A, \star, \alpha)$ with nondegenerate multiplication and surjective $\alpha$ is always a twist. Obviously, to properly understand this statement, a definition of the concept of nondegeneracy used is needed:
\begin{definition}
Let $(A, \star)$ be a set together with a binary operation $\star:A \times A \rightarrow A$. Then $(A, \star)$ is called \emph{left-degenerate} if there exist $a \neq b \in A$ such that $x \star a = x \star b$ for all $x \in A$. Right degeneracy is defined accordingly. $A$ is called \emph{two-sided degenerate} if it is both right and left degenerate. It is called \emph{strongly degenerate} if there exist $a \neq b \in A$ such that both $x \star a = x \star b$ and $a \star x = b \star x$ for all $x \in A$.
\end{definition}
If for example $(A, \star, +)$ is a nonassociative ring without unit, left degeneracy in the sense defined above means that there exists some $0 \neq c \in A$ such that $x \star c = 0$ for all $x \in A$.\\
We are now ready to state and prove:
\begin{proposition} \label{twisting-theorem}
Let $(A, \star, \alpha)$ be a hom-associative structure with $\alpha$ surjective. Then either:
\begin{enumerate}
\item $A$ is a twist. 
\item $A$ is strongly degenerate.
\end{enumerate}
\end{proposition}
\begin{proof} Take some $\beta: A \rightarrow A$ such that $\alpha \circ \beta = id_A$ and define $x \cdot y := \beta(x \star y)$ for $x,y \in A$. Assume that $a \star \beta(b \star c) \neq \beta(a \star b) \star c$ for some $a,b,c \in A$, i.e. $(A, \cdot)$ not associative and suppose that $(A, \star)$ is not strongly degenerate. Then setting $r := a \star \beta(b \star c), s := \beta(a \star b) \star c$ we can find some $b \in A$ such that $b \star r \neq b \star s$ or $r \star b \neq s \star b$. Assume without loss of generality the former. Then by repeating the same argument, we find $a \in A$ such that either $a \star (b \star r) \neq a \star (b \star s)$ or $(b \star r) \star a \neq (b \star s) \star a$. But both inequalities are in contradiction to Lemma \ref{haupt-hilfssatz}, so the proposition follows.
\end{proof}
Using similar ideas, it is possible to show also the following observation on properties of the twisting map:
\begin{proposition} \label{injectivity-theorem}
Suppose $(A, \star, \alpha)$ hom-associative, not strongly degenerate, and $\alpha$ surjective. Then $\alpha$ is in fact bijective.
\end{proposition}
\begin{proof}
Define $\beta$ as in the proof of the previous proposition and suppose that there is $\xi \in A$ with $\beta(\alpha(\xi)) \neq \xi$. As before, we can assume without loss of generality that there exists $b \in A$ with $b \star  \beta(\alpha(\xi)) \neq b \star \xi$ and can then find an $a \in A$ with either $a \star (b \star \beta(\alpha(\xi))) \neq a \star (b \star \xi)$ or $(b \star \beta(\alpha(\xi))) \star a \neq (b \star \xi) \star a$. The first of these possibilities leads to a contradiction due to the general identity
\begin{equation}
x \star (y \star \beta(\alpha(z))) = (\beta(x) \star y) \star \alpha(z) = x \star (y \star z). \label{eqntmp}
\end{equation}
The second case requires application of the same line of reasoning again. We can find $c \in A$ such that either
\begin{equation}
c \star ((b \star \beta(\alpha(\xi))) \star a) \neq c \star ((b \star \xi) \star a) \label{case-1}
\end{equation}
or
\begin{equation}
((b \star \beta(\alpha(\xi))) \star a) \star c \neq ((b \star \xi) \star a) \star c.\label{case-2}
\end{equation}
We will show that both of these possibilities are in contradiction to general identities on $A$. As far as (\ref{case-1}) is concerned, we find that
\begin{eqnarray*}
c \star ((b \star \xi) \star a) & = & (\beta(c) \star (b \star \xi)) \star \alpha(a)\\
& \overset{\text{Eq. \ref{eqntmp}}}{=} & (\beta(c) \star (b \star \beta(\alpha(\xi)))) \star \alpha(a) \\
& = & c \star (( b \star \beta(\alpha(\xi))) \star a).
\end{eqnarray*}
To dispose of (\ref{case-2}), we calculate
\begin{eqnarray*}
((b \star \xi) \star a) \star c & = & (\alpha(b) \star (\xi \star \beta(a))) \star c\\
& \overset{\text{Eq. \ref{eqntmp}}}{=} & (\beta(\alpha^2(b)) \star (\xi \star \beta(a))) \star c\\
& = & \alpha^2(b) \star ((\xi \star \beta(a)) \star \beta(c)) \\
& = & \alpha^2(b) \star ((\xi \star \beta(a)) \star \alpha(\beta^2(c))) \\
& = & \alpha^2(b) \star (\alpha(\xi) \star (\beta(a) \star \beta^2(c)))\\
& = & \alpha^2(b) \star ((\beta(\alpha(\xi)) \star \beta(a)) \star \beta(c))\\
& = & (\alpha(b) \star (\beta(\alpha(\xi)) \star \beta(a))) \star c\\
& = & ((b \star \beta(\alpha(\xi))) \star a) \star c.
\end{eqnarray*}
So both (\ref{case-1}) and (\ref{case-2}) lead to a contradiction. This concludes the proof.
\end{proof}
It is clear that in general the assumption of surjectivity can not be weakened. For instance, $\mathbb{N}$ with addition as binary operation and $\alpha(x) := x + 1$ is clearly hom-associative, but does not arise by the construction described from anything else. The reason is that $0 = 0+0$ is outside the image of $\alpha$, so addition can not be written as an $\alpha$-twisted version of any other operation. Also non-associative examples of this situation can be constructed.\\
If $(A, \star, \alpha, c)$ is left weakly unital, and if $\alpha$ is bijective, we know that $(A, \star)$ is nondegenerate since $c$ acts bijectively by left multiplication. Hence, in this case, Proposition \ref{twisting-theorem} recovers the result from \cite{FG2009} that $A$ is a twisting of an associative structure.
\section{Hom-associative deformation theory} 
We will now explore some applications of our results in the previous section to the deformation theory of Hom-associative algebras. The basic idea we will follow is that both bijectivity of a twisting map and nondegeneracy of a multiplication are properties which are preserved under formal deformation. This enables us to partially ``pull back'' the deformation problem for hom-associative algebras to the deformation problem for associative algebras, which is much better understood.\\
The notion of formal deformations of associative algebras goes back to \cite{Gerstenhaber}. It was extended to hom-algebras in \cite{MS2007}. It is well-known that infinitesimal deformations and obstruction theory of associative deformations are controlled by second and third Hochschild cohomology respectively. Equivalence classes of hom-associative deformations of hom-associative algebras have similarly been identified with elements of a second cohomology module \cite{MS2007}, but so far no cohomology theory for hom-associative algebras has been constructed that would allow a cohomological description of obstruction theory.\\
Throughout this section, $k$ is a commutative ring and $(A, \star, \alpha)$ is a hom-associative $k$-algebra, unless explicitly stated otherwise with nondegenerate multiplication and surjective $\alpha$. By a ``nondegenerate'' multiplication, we will in the sequel always mean a \emph{not strongly degenerate} one.\\
This section is divided into two subsections. The first one briefly recalls the notion of hom-associative formal deformation as given in \cite{MS2007}. In the second subsection, we show that hom-associative deformation preserves non-degeneracy and surjectivity of the twisting map. We use this fact to deduce that deformations of $(A,\star,\alpha)$ have an associative untwist. We prove that this untwist is, in turn, an associative formal deformation of the untwist of $A$.
\subsection{Hom-associative formal deformations}
Let $(A,\star, \alpha)$ be an arbitrary hom-associative algebra. Then \cite{MS2007} give the following definition of a \emph{formal deformation} of $A$:
\begin{definition} \label{formal deformation}
Let $A[[t]]$ be the module of formal power series over $A$ in one variable. Consider a $k[[t]]$-bilinear extension $\mu_t$ of a $k$-bilinear map of type $A \otimes A \rightarrow A[[t]]$ of the form
\begin{displaymath}
\mu_t = \sum_{i \geq 0} t^i \mu_i
\end{displaymath}
with $\mu_0(a,b) = a \star b$ for all $a,b \in A$ and $\mu_i:A \otimes A \rightarrow A$ a bilinear map for every $i \in \mathbb{N}$. Suppose further that we have given a $k[[t]]$-linear map $\alpha_t$ arising by $k[[t]]$-linear extension of a $k$-linear map of the form $\sum_{i \geq 0} X^i \alpha_i$, with $\alpha_0 = \alpha$. Then $(A[[t]], \mu_t, \alpha_t)$ is called a \emph{formal deformation} of $(A, \star, \alpha)$ if $(A[[t]], \mu_t, \alpha_t)$ is hom-associative. 
\end{definition}
In \cite{MS2007}, also a notion of formal equivalence for deformations of hom-associative algebras is defined:
\begin{definition}
Suppose that $(A[[t]], \mu_t, \alpha_t)$ and $(A[[t]], \mu_t', \alpha_t')$ are hom-associative deformations of the hom-associative algebra $(A, \mu, \alpha)$. Then both deformations are called \emph{equivalent} if there exists a formal isomorphism between them, i.e. a $k[[t]]$-linear map $\varphi_t$, compatible with both the deformed multiplications and the deformed twisting maps, of the form
\begin{displaymath}
\varphi_t = \sum_{i \geq 0} t^i \varphi_i
\end{displaymath}
where the $\varphi_i$ are linear maps $A \rightarrow A$ and $\varphi_0 = id_A$. Compatibility with the deformed multiplications means that $\varphi_t \circ \mu_t = \mu_t' \circ (\varphi_t \otimes \varphi_t)$, compatibility to the twisting maps means $\varphi_t \circ \alpha_t = \alpha_t' \circ \varphi_t$.
\end{definition}
\subsection{The nondegenerate, surjective twisting case}
Now if $(A, \star, \alpha)$ is a nondegenerate hom-associative algebra with surjective twisting, we know by the results of the first section that $\alpha$ is in fact a bijection. Consider then a hom-associative deformation $(A, \mu_t, \alpha_t)$ of $(A, \star, \alpha)$. Since in
\begin{displaymath}
\alpha_t = \sum_{i \geq 0} t^i \alpha_i
\end{displaymath}
we have $\alpha_0 = \alpha$ by definition, the usual arguments on invertibility of formal power series yield bijectivity of $\alpha_t$ immediately.\\
Nondegeneracy of the multiplication is also preserved under hom-associative deformation. To see this, let $a := \sum_{i \geq 0} t^i a_i$ be a nonzero element of $A[[t]]$. Choose $n \in \mathbb{N}$ such that $a_i = 0$ for all $i < n$ and $a_n \neq 0$. Denote by $\star_t$ a formal deformation of the original product. Since $A$ was nondegenerate, we can find a $b \in A$ such that $a_n \star b \neq 0$ or $b \star a_n \neq 0$. Assume without loss of generality $b \star a_n \neq 0$. Then since $b \star_t a = t^n b \star a_n + [\text{terms of order $\geq n+1$}]$ we have also $b \star_t a \neq 0$.\\
We have therefore proven:
\begin{quote}
Formal hom-associative deformations of not strongly degenerate hom-associative algebras with surjective twisting map are twists.
\end{quote}
Consider now a formal deformation $(A[[t]], \mu_t, \alpha_t)$ of $(A, \mu, \alpha)$ under these conditions. Then the untwist of the deformed algebra has $\alpha_t^{-1} \circ \mu_t$ as multiplication. This can be expressed as a $k[[t]]$-linear extension of a formal power series with order zero term $\alpha^{-1} \circ \mu$, which means that the untwist of the deformed algebra is an associative formal deformation of $(A, \alpha^{-1} \circ \mu)$, the untwist of $A$. Set the following:
\begin{definition}
Let $(A, \mu)$ be an associative algebra, let $(A[[t]], \star_t)$ be an associative formal deformation of $A$ and let $\alpha_t:A[[t]] \rightarrow A[[t]]$ be a $k[[t]]$-linear map of the form
\begin{displaymath}
\alpha_t = \sum_{i \geq 0} t^i \alpha_i,
\end{displaymath}
with the $\alpha_i$ being module endomorphisms of $A$ and such that
\begin{equation} \label{formal twist condition}
\alpha_t(\alpha_t(x) \star_t \alpha_t(y \star_t z)) = \alpha_t(\alpha_t(x \star_t y) \star_t \alpha_t(z)). 
\end{equation}
Then $\alpha_t$ is called a \emph{formal twisting} of $A[[t]]$ with respect to $\mu_t$.
\end{definition}
Three remarks about Eq. \ref{formal twist condition} are in order. First, Eq. \ref{formal twist condition} is obviously designed in such a way as to give rise to a hom-associative twisting of the formal deformation $(A, \star_t)$. Second, using the $k[[t]]$-linearity of all maps appearing in a standard way, one can check that it is sufficient to verify Eq. \ref{formal twist condition} for $x,y,z \in A$. Third, if $(A, \mu, \alpha)$ is hom-associative and nondegenerate with surjective twisting, then our previous observation that any hom-associative deformations of $A$ can be obtained as twists of associative deformations of $(A, \alpha^{-1} \circ \mu)$ is immediately refined to
\begin{proposition} \label{first deformation proposition}
Any formal hom-associative deformation $(A[[t]], \mu_t, \alpha_t)$ of $(A, \mu, \alpha)$ is obtained from a formal twisting with degree zero component $\alpha$ of an associative deformation of $(A, \alpha^{-1} \circ \mu)$.
\end{proposition}
These findings suggest treating the deformation problem for a nondegenerate hom-associative algebra $(A, \star, \alpha)$ with surjective $\alpha$ in the following way:
\begin{enumerate}
\item Compute $(A, \alpha^{-1} \circ \mu)$.
\item Use associative deformation theory to classify deformations of this algebra.
\item Finally, find all formal twistings with degree zero component $\alpha$ of these deformations.
\end{enumerate}
There are two problems standing in the way of this program:
\begin{enumerate}
\item One needs to verify that any formal twisting with degree zero component $\alpha$ of an associative formal deformation of $(A, \alpha^{-1} \circ \mu)$ gives rise to a hom-associative formal deformation of $(A, \mu, \alpha)$. This is not hard, since hom-associativity of any such formal twisting is true by construction and because verification of the rest of the ``formal deformation'' condition involves only calculations in degree zero terms.
\item One should check that, in order to find \emph{all} hom-associative formal deformations of the original algebra, it suffices to carry out the last step for \emph{one} member of each equivalence class of formal deformations of $(A, \alpha^{-1} \circ \mu)$.
\end{enumerate}
We will deal with the second problem now. We start with the following:
\begin{remark}
Let $(A, \cdot, \alpha)$ be an associative algebra together with a $k$-module homomorphism $\alpha$ satisfying $\alpha(\alpha(x)\alpha(yz)) = \alpha(\alpha(xy)\alpha(z))$. Assume that $(A, \cdot')$ is another $k$-algebra structure isomorphic to $(A, \cdot)$ via $\varphi:A \rightarrow A$. Then $\alpha' := \varphi \circ \alpha \circ \varphi^{-1}$ is a twisting for $(A, \cdot')$ and the hom-associative algebras induced by $\alpha$ on $(A, \cdot)$ and by $\alpha'$ on $(A, \cdot')$ are isomorphic as hom-algebras.
\end{remark}
\begin{proof}
It is clear that $\varphi \circ \alpha = \alpha' \circ \varphi$. Next, we need to prove that $\alpha'$ actually induces a twisted multiplication on $(A, \cdot')$ which with respect to $\alpha'$ is hom-associative. To do this, we calculate
\begin{eqnarray*}
\alpha'(\alpha'(x) \cdot' \alpha'(y \cdot' z)) & = & \alpha'(\varphi(\alpha(\varphi^{-1}(x))) \cdot' \varphi(\alpha(\varphi^{-1}(y \cdot' z))))\\
& = & \alpha'(\varphi(\alpha(\varphi^{-1}(x)) \cdot \alpha(\varphi^{-1}(y) \cdot \varphi^{-1}(z)))) \\
& = & \varphi(\alpha(\alpha(\varphi^{-1}(x)) \cdot \alpha(\varphi^{-1}(y) \cdot \varphi^{-1}(z))))\\
& = & \varphi(\alpha(\alpha(\varphi^{-1}(x) \cdot \varphi^{-1}(y)) \cdot \alpha(\varphi^{-1}(z))))\\
& = & \varphi(\alpha(\alpha(\varphi^{-1}(x \cdot' y)))) \cdot' \varphi(\alpha(\varphi^{-1}(z)))\\
& = & \alpha'(\varphi(\alpha(\varphi^{-1}(x \cdot' y)))) \cdot' \alpha'(z)\\
& = & \alpha'(\alpha'(x \cdot' y) \cdot' \alpha'(z)).
\end{eqnarray*}
What remains to be shown is compatibility of the isomorphism $\varphi$ with the hom-associative multiplications $x \star y := \alpha(x \cdot y)$ and $x \star' y := \alpha'(x \cdot' y)$. This is done by calculating
\begin{displaymath}
\varphi(x \star y) = \varphi(\alpha(x \cdot y)) = \alpha'(\varphi(x \cdot y)) = \alpha'(\varphi(x) \cdot' \varphi(y)) = \varphi(x) \star' \varphi(y).
\end{displaymath}
\end{proof}
It is clear that the previous remark holds also when we do everything on the formal level, i.e. replace isomorphisms with formal isomorphisms and twistings with formal twistings. Only one thing still needs to be checked. Assume that $(A, \cdot, \alpha)$ is an associative algebra together with a twisting $\alpha$ satisfying $\alpha(\alpha(x) \cdot \alpha(y \cdot z)) = \alpha(\alpha(x \cdot y) \cdot \alpha(z))$. Suppose further that $(A[[t]], \star)$ and $(A[[t]], \star')$ are associative formal deformations of $A$, that $\varphi_t$ is a formal isomorphism between them and that $\alpha_t$ is a deformation compatible with $(A[[t]], \star)$ of the twisting $\alpha$. Then for our deformation program to work, we must verify that $\varphi_t \circ \alpha_t \circ \varphi_t^{-1}$ has $\alpha$ as degree zero contribution. But this follows from the fact that $\varphi$ is by definition a deformation of the identity map.
\subsection{Acknowledgements}
The author is grateful to Yael Fregier for useful discussions and for proofreading of parts of an earlier version of this paper. He also wishes to thank the University of Luxembourg for providing excellent working conditions.\\
In the research leading to this work, we were aided by the computer programs Prover9 and Mace4 written by William McCune \cite{McCune}. Specifically, the proofs of Propositions \ref{twisting-theorem} and \ref{injectivity-theorem} were obtained by generalisation, structuring and (minor) simplification of proof objects provided by Prover9 for the case of a hom-associative structure $A$ containing some $c \in A$ such that $c \star x = c \star y$ only if $x = y$. Prover9 is also capable of proving suitable formalisations of these propositions directly, but the proofs given in this paper are in our view much easier to understand than the ones obtained in this way. Mace4 was very helpful in producing finite examples of hom-associative structures which guided our search for general structure theorems.
\bibliography{twist}
\end{document}